\documentstyle[12pt]{amsart}
\newcommand{\catqot}{/\hskip-3pt/}
\newcommand{\C}{{\Bbb C}}
\newcommand{\E}{{\cal E}}
\newcommand{\End}{\mathop{\rm End}}
\newcommand{\F}{{\cal F}}

\newcommand{\id}{\mathop{\rm id}}

\renewcommand{\L}{{\cal L}}
\newcommand{\la}{\lambda}
\newcommand{\M}{{\cal M}}
\renewcommand{\O}{{\cal O}}
\renewcommand{\P}{{\Bbb P}}
\newcommand{\Pic}{\mathop{\rm Pic}}
\newcommand{\Q}{\Bbb Q}
\newcommand{\SL}{\mathop{\rm SL}}
\newcommand{\Spec}{\mathop{\rm Spec}}
\renewcommand{\tilde}{\widetilde}

\newcommand{\N}{{\Bbb N}}
\newcommand{\lra}{\longrightarrow}

\newcommand{\p}{\prime}
\newcommand{\q}{\quad}
\renewcommand{\phi}{\varphi}
\newcommand{\rk}{\mathop{\rm rk}}
\newcommand{\eps}{\varepsilon}
\newcommand{\Tors}{\mathop{\rm Tors}}
\theoremstyle{plain}
\newtheorem{Thm}{Theorem}[section]
\newtheorem{Cor}[Thm]{Corollary}
\newtheorem{Prop}[Thm]{Proposition}
\newtheorem{Lem}[Thm]{Lemma}

\theoremstyle{remark}
\newtheorem{Rem}[Thm]{Remark}

\title{Projective moduli for Hitchin pairs}
\author{A.\ Schmitt}
\address
{Institut f\"ur Mathematik\\ Universit\"at Z\"urich\\
Winterthurerstr.\ 190\\ CH-8057 Z\"urich\\ Switzerland}
\email{schmitt@@math.unizh.ch}
\begin{document}
\maketitle
\centerline{$\underline{\hskip 5truecm}$}
\vskip 1.2truecm
\section*{Introduction}
In the paper \cite{Hi}, Hitchin studied pairs $(E,\phi)$,
where $E$ is a vector bundle of rank two
with a fixed determinant on a curve $C$  and $\phi\colon E\lra E\otimes K_C$
is a trace free homomorphism, and constructed a moduli space
for them.
This moduli space carries the structure of a non-complete, quasi-projective
algebraic variety.
Later, Nitsure \cite{Ni} gave an algebraic construction of moduli spaces
of pairs $(E,\phi)$ over a curve $C$ 
consisting of a vector bundle $E$ of fixed degree and rank
and a homomorphism $\phi\colon E\lra E\otimes L$ where $L$ is some previously
chosen line bundle.
He also obtained non-complete moduli spaces.
The most general results were obtained by Yokogawa \cite{Yo}.
In his paper, $C$ is replaced by a relative scheme $f\colon X\lra S$ where
$f$ is a smooth, projective, geometrically integral morphism
and $S$ is a scheme of finite type over a universally Japanese ring,
and $L$ by a locally free sheaf $F$ on $X$.
\par
It is the aim of our paper to compactify some of the spaces
obtained by Yokogawa, namely those where $S=\Spec\C$ and
$F$ is again a line bundle.
In order to avoid confusion with the objects studied e.g.\ by Simpson,
we will call our objects \it (oriented) Hitchin pairs\rm .
\par
We shall also mention that,
only recently, T.\ Hausel compactified the space of oriented Hitchin pairs 
of rank two
with fixed determinant over a curve $C$, using methods from symplectic
geometry. This result and a detailed investigation of the resulting spaces
will appear in a forthcoming preprint of his.
\par
The structure of this note is as follows:
In the first section we treat the case where $X$ is a point.
This case shows how to define Hitchin pairs correctly and suggests the
definition of (semi)stability.
Then we prove a boundedness result following \cite{Ni}, 
construct a projective parameter space for semistable Hitchin pairs
and a universal family on this parameter space,
and finally define a linearized $\SL(V)$-action on this parameter 
space such that the moduli space is given as $\text{parameter space}\catqot
\SL(V)$.
After these constructions, we prove the (semi)stability criterion.
\par
At some places, the techniques of our notes are similar to
those in \cite{OST}.
Hence, we often omit or sketch only briefly arguments which were carried
out in detail in \cite{OST} in an analogous situation.
\section*{Acknowledgements}
I want to thank Professor Okonek and Dr.\ A.\ Teleman for suggesting
the problem and discussing various details of the proof with me.
The author acknowledges support by AGE --- Algebraic Geometry in 
Europe,
contract No.\ ER-BCHRXCT 940557 (BBW 93.0187), 
and by SNF, Nr.\ 2000 -- 045209.95/1.
\section{Compactifying the categorical quotient of
a vector space}
\label{Simple}
Let $G$ be a reductive algebraic group acting linearly on a vector space $V$.
Consider the categorical quotient $W:=V\catqot G=\Spec\C[V]^G$.
The torus $\C^*$ acts canonically on $V$, and this action
commutes with the given action of $G$.
Now, let $G$ act trivially on $\C$ and let $\C^*$ act on $\C$ by 
multiplication.
We obtain a ($G\times\C^*$)-action on $V\oplus \C$.
Observe that the equivalence relation induced by the given action
is the following:
$$
(v_1,\eps_1)\q\sim\q (v_2,\eps_2)
\qquad\Leftrightarrow\qquad \exists z\in\C^*,\
g\in G: v_2= z\cdot (g\cdot v_1);\ \eps_2=z\cdot\eps_1.
$$
The point is 
$(v,\eps)\in V\oplus \C$ is semistable if and only if 
$[v,\eps]\in\P(V\oplus \C)$ is $G$-semistable.
By the Hilbert criterion, the latter happens if and only if
either $\eps\neq 0$ or $v\in V$ is $G$-semistable.
The space $(V\oplus\C)\catqot (G\times\C^*)=\P(V\oplus \C)\catqot G$
obviously is a projective variety containing $W$ as an open affine subvariety.
Let $W^{ss}$ be the image of the $G$-semistable points in $V$.
The $\C^*$-action on $V$ induces a $\C^*$-action on $W^{ss}$.
We observe that we have compactified $W$ with $W^{ss}\catqot\C^*$.
\par
Applying the above discussion to the case $G=\SL_n(\C)$ and $V=M_n(\C)$
(this is the case of Hitchin pairs over a point) shows that
$(m,\eps)\in M_n(\C)\oplus\C$ is semistable if and only if
either $\eps\neq 0$ or $m$ is not nilpotent.
\begin{Rem}
Comparing this with \cite{Ni}, Thm.2.8, for $r=p$ and $N=1$
shows that the semistability criterion stated there is
false for points at infinity.
\end{Rem}
Our general construction is basically a relative version of the above
over a projective scheme.
\section{Hitchin pairs}
Throughout this paper, we will work over the field of complex numbers.
Let $X$ be a smooth projective variety of dimension $n$. If $n>1$, we fix 
an ample divisor $H$ on $X$ whose associated line bundle will be denoted 
by $\O_X(1)$. We will use $H$ to compute degrees and Hilbert polynomials.
The Hilbert polynomial of a coherent sheaf $\F$ will be denoted by $P_\F$.
We also fix a line bundle $L$ and a Hilbert polynomial $P$.
The degree and the rank given by $P$ will be denoted by $d$ and $r$,
respectively.
Let $\Pic(X)$ be the Picard scheme of $X$. We fix a Poincar\'e sheaf
$\L$ on $\Pic(X)\times X$.
Furthermore, for a coherent sheaf $\E$, set 
$\L[\E]:=\L_{|\{[\det\E]\}\times X}$.
This sheaf depends only on the isomorphy class of $\E$.
Unlike the situation in \cite{OST}, the sheaf $\L$ will play
no essential r\^ole
in our considerations.
\subsection{Oriented  Hitchin Pairs}
\label{or}
An \it oriented Hitchin pair of type $(\L,P,L)$ \rm is a triple 
$(\E,\sigma,\phi)$
consisting of a torsion free coherent sheaf $\E$ with $P_\E=P$,
a homomorphism $\sigma\colon \det(\E)\lra \L[\E]$, and a homomorphism 
$\phi\colon \E\lra\E\otimes L$.
Two oriented Hitchin pairs $(\E_1,\sigma_1,\phi_1)$
and $(\E_2,\sigma_2,\phi_2)$ of type $(P,\L,L)$
are called \it equivalent\rm , if there 
is an isomorphism $\psi\colon \E_1\lra\E_2$ 
such that
$\phi_2\circ\psi=(\psi\otimes \id_L)\circ \phi_1$
and $\sigma_1=\sigma_2\circ\det\psi$.
\begin{Rem}
Of course, we can fix a line bundle $\L_0$ on $X$ and consider 
oriented Hitchin pairs $(\E,\sigma,\phi)$ such that $\det\E\cong\L_0$.
Our proofs carry over to this situation.
\end{Rem}
Now, consider pairs $(\E,\phi)$ consisting of a torsion free coherent
sheaf
$\E$ with $P_\E=P$ and a homomorphism $\phi\colon \E\lra\E\otimes L$.
We say that $(\E_1,\phi_1)$ is \it equivalent to \rm $(\E_2,\phi_2)$
if and only if there is an isomorphism $\psi\colon \E_1\lra\E_2$
fulfilling $\phi_2\circ\psi= (\psi\otimes\id_L)\circ\phi_1$.
Given a pair $(\E,\phi)$, we can choose a non-zero orientation 
$\sigma\colon\det\E\lra\L[\E]$ in order to obtain an 
oriented Hitchin pair. We observe that the equivalence class
of $(\E,\sigma,\phi)$ does not depend on the choice of the orientation
$\sigma$.
Therefore, we call a pair $(\E,\phi)$ as above \it an oriented
Hitchin pair of type $(L,P)$\rm .
Let $S$ be a noetherian scheme. \it A family of oriented Hitchin pairs
of type $(L,P)$
parametrized by $S$ \rm is a pair $({\frak E}_S,\phi_S)$ where 
${\frak E}_S$ is a coherent sheaf on $S\times X$ and $\phi_S$
is an element of 
$H^0(S\times X, \underline{\End}({\frak E}_S)\otimes \pi_X^*L)$
such that $({\frak E}_{S|\{s\}\times X},\phi_{S|\{s\}\times X})$
is an oriented Hitchin pair of type $(L,P)$ for any closed point $s\in S$.
Two families $({\frak E}^i_S,\phi^i_S)$, $i=1,2$, are said to be
\it equivalent\rm , if there is an isomorphism $\psi_S\colon {\frak E}^1_S
\lra {\frak E}_S^2\otimes\pi_X^*L$ with 
$\phi_S^2\circ\psi_S=\bigl((\psi_S\otimes\id_{\pi_X^*L})\circ \phi_S^1\bigr)$.
\subsection{Hitchin Pairs}
A \it Hitchin pair of type $(L,P)$ \rm is a triple $(\E,\eps,\phi)$
consisting of a torsion free coherent sheaf $\E$ with $P_\E=P$,
a complex number $\eps\in\C$, and a homomorphism 
$\phi\colon \E\lra\E\otimes L$.
Two Hitchin pairs $(\E_1,\eps_1,\phi_1)$
and $(\E_2,\eps_2,\phi_2)$ are called \it equivalent\rm , if there 
are an isomorphism $\psi\colon \E_1\lra\E_2$ and a complex number
$z\in\C^*$ such that 
$\phi_2\circ\psi=\bigl(\psi\otimes (z\cdot \id_L)\bigr)\circ \phi_1$
and $\eps_2=z\eps_1$.
Let $S$ be a noetherian scheme.
\it A family of Hitchin pairs of type $(L,P)$ parametrized by $S$ \rm
is a quadruple $({\frak E}_S,\eps_S,\phi_S,{\frak M}_S)$
consisting of a coherent sheaf ${\frak E}_S$ over $S\times X$,
an invertible sheaf ${\frak M}_S$ over $S$, 
a section $\eps_S\in H^0(S,{\frak M}_S)$, and an element
$\phi_S\in H^0(S\times X, 
\underline{\End}({\frak E}_S)\otimes\pi_S^*{\frak M}_S\otimes \pi_X^*L)$
such that its restriction to $\{s\}\times X$ is a Hitchin pair
of type $(L,P)$ for any closed point $s\in S$.
The family $({\frak E}^1_S,\eps^1_S,\phi^1_S,{\frak M}^1_S)$
is said to be \it equivalent to \rm the family
$({\frak E}^2_S,\eps^2_S,\phi^2_S,{\frak M}^2_S)$
if there are isomorphisms $\psi_S\colon {\frak E}_S^1\lra {\frak E}_S^2$
and $z_S\colon {\frak M}_S^1\lra {\frak M}_S^2$ with
$\phi_S^2\circ \psi_S=\bigl(\psi_S\otimes \pi_S^*z_S\otimes 
\id_{\pi_X^*L}\bigr)\circ\phi_S^1$ and $\eps_S^2=\eps_S^1\circ z_S$.
\subsection{(Semi)Stability}
We call a Hitchin pair $(\E,\eps,\phi)$ of type $(L,P)$ \it (semi)stable\rm , 
if the
following two conditions are satisfied:
\begin{enumerate}
\item For any $\phi$-invariant subsheaf $0\neq \F\subset\E$ we have:
$(P_\F/\rk\F)\quad (\le)\quad (P/r)$.
\item Either $\eps\neq 0$, or 
$(\phi\otimes\id_{L^{\otimes r-1}})\circ\cdots\circ \phi\neq 0$.
\end{enumerate}
\begin{Rem}
As usual, there are the corresponding notions of \it slope-(semi)stability\rm .
Slope-stability implies stability and semistability implies
slope-semistability.
\end{Rem}
We are now able to define the functors $M^{(s)s}_{(L,P)}$ of equivalence 
classes of families of (semi)stable Hitchin pairs of type $(L,P)$.
The functors of families of (semi)stable oriented Hitchin pairs of type
$(L,P)$
are the open subfunctors $\eps\neq 0$.
\section{Boundedness}
It is the aim of this section to show that the family of isomorphy classes
of torsion free coherent sheaves occuring in slope-semistable
Hitchin pairs of type $(L,P)$ is bounded.
We recall that any torsion free coherent sheaf $\E$ possesses
a Harder-Narasimhan filtration
$$0=\E_0\subset\E_1\subset\cdots\subset\E_l=\E$$
where $\E_i/\E_{i-1}$ is the subsheaf of maximal rank of $\E/\E_{i-1}$
for which $P_{\E_i/\E_{i-1}}$ is maximal.
We have
\begin{equation}
\label{G1}
\mu(\E_i/\E_{i-1})\qquad \ge \qquad \mu(\E_{i+1}/\E_i),\qquad i=1,...,l-1.
\end{equation}
A simple inductive argument shows that $\mu(\E_i)\ >\ \mu(\E)$, $i=1,...,l$, 
when
$\mu(\E_1)\ >\ \mu(\E)$.
By a theorem of Maruyama \cite{Ma}, it is enough to bound $\mu(\E_1)$
for torsion free coherent sheaves occuring in slope-semistable
Hitchin pairs of type $(L,P)$:
\begin{Thm}
\label{B1}
For any torsion free coherent sheaf which is part of a slope-semistable
Hitchin pair of type $(L,P)$, we have
$$\mu(\E_1)\qquad \le\qquad \max\Bigl\{\, \mu(\E),\ \mu(\E)+{(r-1)^2\over r}\deg L\,\Bigr\}.$$
\end{Thm}
\begin{pf}
We follow the proof of \cite{Ni}, Prop.3.2, in the case of curves.
Let $(\E,\eps,\phi)$ be a slope-semistable Hitchin pair of type $(L,P)$.
If $\mu(\E_1)\ \le\ \mu(\E)$, there is nothing to show.
Otherwise, as we have seen above, $\mu(\E_i)\ >\ \mu(\E)$, $i=1,...,l$.
By definition of slope-semistability, this means that the $\E_i$ are not
$\phi$-invariant.
Hence, the homomorphism $\phi_i\colon \E_i\lra \E/\E_i\otimes L$
is not trivial for $i=1,...,l$.
Let $\iota\in\{\, 0,...,i-1\,\}$ be maximal with $\phi_i(\E_{\iota})=0$
and $\kappa\in \{\, i+1,...,l\,\}$ minimal with 
$\phi_i(\E_i)\subset (\E_\kappa/\E_i)\otimes L$.
With these choices, the induced homomorphism from $\E_{\iota+1}/\E_\iota$
to $\E_\kappa/\E_{\kappa-1}\otimes L$ is non-trivial.
Both of these sheaves are slope-semistable, so that $
\mu(\E_{\iota+1}/\E_\iota)
\ \le\ \mu(\E_\kappa/\E_{\kappa-1}) + \deg L$.
By (\ref{G1}), we get
\begin{equation}
\label{G2}
\mu(\E_i/\E_{i-1})\qquad \le\qquad \mu(\E_{i+1}/\E_i)+\deg L.
\end{equation}
Summing these inequalities from $i=1$ to $i=l-1$ yields
$$\mu(\E_1)\q \le\q \mu(\E/\E_{l-1})+(l-1)\deg L\q \le\q \mu(\E/\E_{l-1})+(r-1)\deg L.$$
Since $\mu(\E_1)+(r-1)\mu(\E/\E_{l-1})\ \le\ r\mu(\E)$, i.e.,
$$
\mu(\E/\E_{l-1})\qquad \le \qquad {d-\mu(\E_1)\over r-1},
$$
the assertion of the theorem follows.
\end{pf}
\begin{Rem}
\label{B2}
Fix a number $m$ such that $L\subset\O_X(m)$.
For any coherent sheaf $\F$, we obviously have $P_{\F\otimes L}\ \le\
P_{\F(m)}$.
It is easy to see that there is a constant $C^\p$ \sl depending
only on $H$ and $m$ \rm with $P_{\F(m)}\ \le\ P_\F+C x^{n-1}$.
\end{Rem}
We can now carry out the proof of \ref{B1} for semistable Hitchin pairs
and Hilbert polynomials, where we replace (\ref{G2}) by
$$
{P_{\E_i/\E_{i-1}}\over \rk\E_i-\rk\E_{i-1}}\qquad \le\qquad
{P_{\E_{i+1}/\E_i}+ C x^{n-1}\over \rk\E_i-\rk\E_{i-1}}
\qquad\le\quad {P_{\E_{i+1}/\E_i}\over \rk\E_i-\rk\E_{i-1}}+C x^{n-1}.
$$
This gives
$$
{P_{\E_1}\over \rk\E_1}\qquad\le\qquad {P +(r-1)^2 C x^{n-1}\over r}.
$$

\section{A parameter space for semistable Hitchin pairs}
For $\mu\in\N$, we define $P_\mu$ by $P_\mu(x):=P(x+\mu)$.
Twisting by $\O_X(\mu)$ yields an isomorphism between the functors
$M^{(s)s}_{(L,P)}$ and $M^{(s)s}_{(L,P_\mu)}$.
By Theorem~\ref{B1}, we may assume that any torsion free coherent sheaf
$\E$ appearing in a semistable Hitchin pair of type $(L,P)$ fulfills the 
following conditions:
\begin{enumerate}
\item $\E$ is globally generated.
\item $H^i(X,\E)=0$ for every $i>0$.
\end{enumerate}
Let $p:=P(0)$, $V$ be a complex vector space of dimension $p$,
and ${\frak Q}$ the projective 
Quot scheme of (all) quotients of $V\otimes\O_X$
with Hilbert polynomial $P$.
On the product ${\frak Q}\times X$, there is a universal quotient
$$q_{\frak Q}\colon V\otimes \O_{{\frak Q}\times X}\lra {\frak E}_{\frak Q}.$$
We choose $m$ large enough, so that $L\subset \O_X(m)$ and so that
$\O_X(m)$ is globally generated.
Furthermore, we choose $\nu$ large enough, so that $q_{\frak Q}(\nu)$
induces a closed embedding ${\frak Q}\subset {\frak G}:=
\mathop{\rm Gr}\bigl(V\otimes H^0(\O_X(\nu)), P(\nu)\bigr)$
and so that the multiplication map $H^0(\O_X(\nu))\otimes H^0(\O_X(m))\lra
H^0(\O_X(\nu m))$ is surjective.
We set $N:=H^0(\O_X(\nu))$, $M:=H^0(\O_X(m))$, and $W:=V\otimes N$.
By our choice of $\nu$, for any Hitchin pair $(\E,\eps,\phi)$,
the map $V\otimes N\lra H^0(\E(\nu))$ is surjective.
It follows that $\phi\otimes \id_{\O_X(\nu)}$ is induced by
an element $f\in W^\vee\otimes W\otimes M$.
Set ${\bf P}:=\P(\C\oplus W^\vee\otimes W\otimes M^\vee)$, and let
$${\bf s}\colon \O_{\bf P}\lra [\C\oplus W^\vee\otimes W\otimes M]\otimes
\O_{\bf P}(1)$$
be the tautological section.
First, we can construct a subscheme $\tilde{\frak P}\subset {\frak Q}\times
{\bf P}$ whose closed points are those $s=([q],\tilde{s})\in {\frak Q}\times
{\bf P}$ for which the second component of $\pi_{\bf P}^*{\bf s}$
induces a homomorphism ${\frak E}_{{\frak Q}|\{[q]\}\times X}(\nu)
\lra {\frak E}_{{\frak Q}|\{[q]\}\times X}(\nu)\otimes H^m$.
Let ${\frak E}_{\tilde{\frak P}}$ be the restriction of 
$\pi_{\frak Q}^*{\frak E}_{\frak Q}$ to $\tilde{\frak P}\times X$
and
$$
h_{\tilde{\frak P}}\colon {\frak E}_{\tilde{\frak P}}\lra
{\frak E}_{\tilde{\frak P}}\otimes \pi_X^*(\O_X(m)/L)
$$
be the induced homomorphism.
We then define ${\frak P}$ as the closed subscheme of $\tilde{\frak P}$
whose closed points are those $s\in \tilde{\frak P}$ for which
$h_{|\{s\}\times X}\equiv 0$.
The scheme ${\frak P}$ is a parameter space for pairs $([q\colon V\otimes\O_X
\lra\E],[\eps,\phi])$ with $[q]\in {\frak Q}$, $[\eps,\phi]\in \P(\C\oplus
H^0(\underline{\End}\E\otimes L)^\vee)$.
On ${\frak P}\times X$, there exists a universal family 
$({\frak E}_{\frak P},\eps_{\frak P}, \phi_{\frak P}, {\frak M}_{\frak P})$.
\par
Denote by ${\frak P}^{\mathop{\rm iso}}$ the open set of pairs
$([q\colon V\otimes\O_X\lra \E],[\eps,\phi])$ for which $H^0(q)$ is
an isomorphism.
It is not hard to see that any family of semistable Hitchin pairs
of type $(L,P)$ is locally induced by  morphisms to 
${\frak P}^{\mathop{\rm iso}}$.
\section{The $\SL(V)$-action on ${\frak P}$}
On the Quot scheme ${\frak Q}$, there is a natural action
$\rho\colon {\frak Q}\times \SL(V)\lra {\frak Q}$.
Furthermore, there is a natural action of $\SL(V)$
from the right on the vector space $W^\vee\otimes W\otimes M$.
If we let $\SL(V)$ act trivially on $\C$, we get
an action of $\SL(V)$ from the right on the scheme
${\frak Q}\times {\bf P}$.
Finally, we remark that the $\SL(V)$-action leaves the parameter space
${\frak P}$ invariant. Hence, there is an action from the right of
$\SL(V)$ on ${\frak P}$.
We deduce
\begin{Prop}
\label{Op1}
Let $S$ be a noetherian scheme and 
$\beta_i\colon S\lra {\frak P}^{\mathop{\rm iso}}$ two morphisms.
Suppose that the pullbacks via the maps $(\beta_i\times\id_X)$
of the universal family $({\frak E}_{\frak P},\eps_{\frak P},\phi_{\frak P},
{\frak M}_{\frak P})$ are equivalent.
Then there exist an \'etale covering $\tau\colon T\lra S$
and
a morphism $g\colon T\lra\SL(V)$ such that 
$\beta_1\circ \tau=(\beta_2\circ\tau)\cdot g$.
\end{Prop}
\section{The (semi)stable points in ${\frak P}$}
%
%\subsection{The Hilbert-Mumford criterion}
%
Suppose we are given a 
projective scheme $S$ and an action of an algebraic group $G$,
linearized in an invertible sheaf ${\frak M}$.
For a point $s\in S$ and a one parameter subgroup $\la\colon \C^*\lra G$, 
set $s_\infty:=\lim_{z\lra\infty} \la(z)\cdot s$.
Then $s_\infty$ is a fixed point of the $\C^*$-action given by $\la$, 
and $\C^*$ acts on
${\frak M}\otimes\C(s_\infty)$ with weight, say, $\gamma$.
We set $\mu(s,\la):=-\gamma$.
If $G$ is reductive and ${\frak M}$
is ample, then the Hilbert-Mumford criterion says that $s$ is
(semi)stable if and only if $\mu(s,\la)\ (\ge)\ 0$
for every one parameter subgroup $\la$ of $G$.
We will apply this criterion in our situation.
\par
A one parameter subgroup of $\SL(V)$ is determined by the following
data:
\begin{enumerate}
\item A basis $v_1,...,v_p$ of $V$.
\item Weights $\gamma_1\le\cdots\le\gamma_p$ with $\sum_i\gamma_i=0$.
\end{enumerate}
We recall that a weight vector $(\gamma_1,...,\gamma_p)$, satisfying
$\gamma_1\le\cdots\le\gamma_p$ and $\sum_i\gamma_i=0$ is a $\Q$-linear
combination with non-negative coefficients of the weight vectors
$$
\gamma^{(i)}:= (\quad\underbrace{i-p,...,i-p}_{\text{$i$ times}},\underbrace{i,...,i}_{\text{$(p-i)$ times}}\ ).
$$
More precisely,
\begin{equation}
\label{LinKomb}
(\gamma_1,...,\gamma_p)=\sum_{i=1}^{p-1}{\gamma_{i+1}-\gamma_i\over p}\gamma^{(i)}.
\end{equation}
%
%\subsection{(Semi)stability on ${\frak P}$}
%
Let's return to our construction.
Let $\O_{\frak Q}(1)$ be the restriction of the very ample line bundle
on ${\frak G}$ giving the Pluecker embedding.
We denote by $\O(a_1,a_2)$ the restriction of the bundle
$\pi_{\frak Q}^*\O_{\frak Q}(a_1)\otimes \pi_{\bf P}^*\O_{\bf P}(a_2)$
to the parameter space ${\frak P}$.
The $\SL(V)$-action on ${\frak P}$ can be linearized in any of these
sheaves.
We will choose $a_1,a_2>0$ with $a_1<(p-1) a_2$.
For $[q\colon V\otimes \O_X\lra \E]\in {\frak Q}$ and a subspace $U\subset V$,
$\E_U$ is defined to be the subsheaf of $\E$ which is generically generated
by $q(U\otimes\O_X)$ and for which $\E/\E_U$ is torsion free.
Given a basis $v_1,...,v_p$ of $V$, we set 
$\E_i:=\E_{\langle\, v_1,...,v_i\,\rangle}$, so that we obtain a
filtration
$$
\Tors\E=\E_0\subset\E_1\subset\cdots\subset\E_{p-1}\subset\E_p=\E.
$$
Now, either $\E_i=\E_{i+1}$ or $\rk\E_{i+1}=\rk\E_i+1$.
For $\rho=1,...,p$, we set $k_\rho:=\min_{i=1,...,p}\{\,\rk\E_i=\rho\,\}$
and $\underline{k}:=(k_1,...,k_p)$.
Suppose we are given a one parameter subgroup $\la$ of $\SL(V)$.
For a point $[q\colon V\otimes\O_X\lra \E]\in {\frak Q}$,
set $[q_\infty\colon V\otimes\O_X\lra\E]:=\lim_{z\lra\infty}\la(z)\cdot [q]$.
We denote the fibre of $\O_{\frak Q}(a_1)$ over $[q_\infty]$ by $\Lambda$.
Let $v_1,...,v_p$ be a basis of $V$. 
If $\la$ is the one parameter subgroup which is described by the
weight vector $(\gamma_1,...,\gamma_p)$, then $\la$ acts on
$\Lambda$ with weight $a_1\gamma_{\underline{k}}$ (\cite{HL2}, p.309) where
we set
$$
\gamma_{\underline{k}}:=\gamma_{k_1}+\cdots+\gamma_{k_p}.
$$
In particular
$$
\gamma_{\underline{k}}^{(i)}=(p-i)\rk\E_i-i(\rk\E-\rk\E_i)=p\rk\E_i-i\rk\E.
$$
Now, consider the $\SL(V)$-action on ${\bf P}$.
For a point $\tilde{s}\in {\bf P}$ and a one parameter subgroup
$\la$ of $\SL(V)$, define $\tilde{s}_\infty$ as above and let
$E$ be the fibre of $\O_{\bf P}(a_2)$ over $\tilde{s}_\infty$.
For the statement of the next lemma, we need the notion of
a superinvariant subspace which will.
Suppose we are given a homomorphism $f\colon V\otimes N\lra V\otimes N\otimes M$.
A subspace $U\subset V$ is called \it $f$-superinvariant\rm ,
if $U\otimes N\subset\ker f$ and if the induced homomorphism
$\overline{f}\colon (V/U)\otimes N\lra (V/U)\otimes N\otimes M$
is identically zero.
From now on, given an element $s:=([q\colon V\otimes\O_X\lra \E],[\eps,\phi])$
in ${\frak P}$, the associated homomorphism in $W^\vee\otimes W\otimes M$
will be denoted by $f$.
We have the following obvious
\begin{Lem}
\label{SemStab0}
Set $s:=([q\colon V\otimes\O_X\lra \E],[\eps,\phi])$.
The one parameter subgroup $\la$ which is given w.r.t.\ to the basis
$v_1,...,v_p$ by the weight vector $\gamma^{(i)}$ acts on $E$
with weight
\begin{enumerate}
\item $-a_2p$ if $\langle\, v_1,...,v_i\,\rangle$ is not $f$-invariant.
\item $a_2p$ if $\langle\, v_1,...,v_i\,\rangle$ is $f$-superinvariant.
\item $0$ in all the other cases.
\end{enumerate}
\end{Lem}
An immediate consequence is:
\begin{Cor}
\label{SemStab1}
A necessary condition for a point 
$s:=([q\colon V\otimes\O_X\lra \E],[\eps,\phi])$ to be (semi)stable
is that for any $f$-invariant subspace $U\subset V$
$$\dim U\rk\E\q (\le)\q p\rk\E_U.$$
\end{Cor}
\begin{Cor}
\label{SemStab2}
Let $s:=([q\colon V\otimes\O_X\lra \E],[\eps,\phi])$ be a point in ${\frak P}$
and suppose that either $H^0(q)$ is not an isomorphism or that $\E$
is not torsion free.
Then $s$ is not semistable.
\end{Cor}
\begin{pf}
Set $U:=\ker H^0(q)$ in the first case and $U:=H^0(\Tors\E)$ in the second 
case.
Then $U$ clearly violates the condition in Corollary~\ref{SemStab1}.
\end{pf}
We now state the main result of this section:
\begin{Thm}
\label{MainRes}
For $d$ sufficiently large the following assertion holds true:
A point $s:=([q\colon V\otimes\O_X\lra \E],[\eps,\phi])$ is (semi)stable
if and only if $H^0(q)$ is an isomorphism, $\E$ is torsion free,
and $(\E,\eps,\phi)$ is a (semi)stable Hitchin pair.
\end{Thm}
We will need
\begin{Prop}
\label{Aux1}
There is an integer $k_0$ such that for any semistable Hitchin pair,
any subsheaf $\F\subset\E$, and any $k\ge k_0$:
$$
rh^0(\F(k))\q <\q (\rk\F+1) P(k).
$$
\end{Prop}
\begin{pf}
As in the proof on page 305 in \cite{HL2}, we conclude that for any
sufficiently large constant $\kappa$ there is an integer $k_0$
such that for any Hitchin pair $(\E,\eps,\phi)$ and any subsheaf $\F\subset\E$
$$
\text{either}\q |\deg\F-\rk\F\mu(\E)|\ \le\ \kappa\q\text{or}\q h^0(\F(k))/\rk\F\ <\ P(k)/r\q \forall k\ge k_0.
$$
Let ${\frak S}$ be the family of all saturated submodules of torsion free
sheaves $\E$ occuring in the family ${\frak E}_{\frak Q}$ 
which satisfy
$|\deg\F-\rk\F\mu(\E)|\ \le\ \kappa$.
Then this family is bounded (\cite{HL2}, Lemma 2.7).
Hence, we may assume that all $\F\in {\frak S}$ are globally generated and
without higher cohomology.
By the discussions following Remark~\ref{B2}
\begin{eqnarray*}
rh^0(\F(k)) &\le& \rk\F(P(k)+(r-1)^2Ck^{n-1})\\
            &=&  (\rk\F+1)P(k)+[\rk\F(r-1)^2Ck^{n-1}-P(k)]\\
            &\le&  (\rk\F+1)P(k)+[r(r-1)^2Ck^{n-1}-P(k)].
\end{eqnarray*}
Since $C$ does not depend on $d$, we can achieve $[r(r-1)^2Ck^{n-1}-P(k)]<0$
for all $k\ge k_0$.
\end{pf}
We choose $d$ large enough so that $k_0=0$, and so that all modules
$\F$ in the family ${\frak S}$ are globally generated and without 
higher cohomology.
Since there are only finitely many possible Hilbert polynomials
for sheaves in ${\frak S}$, the proof of \ref{Aux1} shows that
we can assume that for any $\F\subset \E$, $\E$ being a torsion free
 member of the family ${\frak E}_{\frak Q}$,
 the inequality $P_\F/\rk\F\ (\le)\ P/r$ is equivalent to the inequality
$h^0(\F)/\rk\F\ (\le)\ p/r$.
\begin{pf*}{Proof of Theorem~\ref{MainRes}}
First, let $([q\colon V\otimes\O_X\lra \E],[\eps,\phi])$ be a (semi)stable
point.
Then, by \ref{SemStab2}, $H^0(q)$ is an isomorphism and $\E$ is torsion free.
Furthermore, \ref{SemStab1} shows that $(\E,\eps,\phi)$
is a (semi)stable Hitchin pair, provided $\eps\neq 0$.
We still have to show that 
$(\phi\otimes\id_{L^{\otimes r-1}})\circ\cdots\circ\phi$ is not zero if
$\eps=0$.
For this, set $\F_i:=\ker[(\phi\otimes\id_{L^{\otimes i-1}})\circ\cdots\circ\phi]$, $i=1,...,r$.
We get a filtration 
$$
0=:\F_0\subset\F_1\subset\cdots\subset\F_{r-1}\subset\F_r=\E
$$
of $\E$.
Choose a basis $v_1,...,v_p$ of $V$ such that there
are $\iota_i$ with $\langle\, v_1,...,v_{\iota_i}\,\rangle=
H^0(\F_i)$, $i=1,...,r$.
Let $\la$ be the one parameter subgroup which is given by the weight vector
$\sum\gamma^{(i)}$.
The assumption $\mu(s,\la)\ (\ge)\ 0$ implies that there is an index $i$
with $a_1(rh^0(\F)-\rk\F p)+a_2 p\ \le\ 0$, in particular
$$
rh^0(\F_i)\qquad <\qquad (\rk\F_i-1)p.
$$
This implies 
$$
rh^0(\E/\F_i)\qquad \ge\qquad (r-\rk\F_i+1)p.
$$
Now, $(\phi\otimes\id_{L^{\otimes i-1}}\circ\cdots\circ\phi)$ maps $\E/\F_i$ 
isomorphically onto a $(\phi\otimes\id_{L^{\otimes i}})$-invariant
subsheaf of $\E\otimes L^{\otimes i}$.
This sheaf can be identified with a 
$(\phi\otimes\id_{H^{\otimes im}})$-invariant
subsheaf of $\E\otimes H^{\otimes im}$.
But $\E\otimes H^{\otimes i m}$ is also semistable, and the assumptions
made before the beginning of the proof hold for this sheaf as well, so that
$$
rh^0(\E/\F_i)\qquad \le\qquad (r-\rk\F_i)P_{\E\otimes H^{\otimes i m}}
=(r-\rk\F_i)P(im),
$$
and, consequently,
$$
(r-\rk\F_i)(P(im)-p)\qquad \le\qquad p.
$$
But when $d$ is large, this is not possible.
\par
Now, we prove the opposite direction: Let $s:=([q\colon V\otimes\O_X\lra \E],
[\eps,\phi])$ be a point such that $H^0(q)$ is an isomorphism and
$(\E,\eps,\phi)$ is a (semi)stable Hitchin pair.
First, suppose $\eps\neq 0$.
Let $v_1,...,v_p$ a basis of $V$.
Let $\lambda$ be given  by the weight vector 
$\gamma=\sum \alpha_i \gamma^{(i)}$.
If all the spaces $\langle\, v_1,...,v_i\,\rangle$ for
which $\alpha_i\neq 0$ are $f$-invariant, then the (semi)stability
condition implies $\gamma_{\underline{k}}\ (\le)\ 0$.
Together with~\ref{SemStab0}, this implies $\mu(s,\la)\ (\ge)\ 0$.
In the other case, let $\alpha$ be the largest coeffictient
of a $\gamma^{(i)}$ for which $\langle\, v_1,...,v_i\,\rangle$
is not $f$-invariant.
By~\ref{Aux1}, $\gamma_{\underline{k}}^{(i)}\ \le\ p$ for $i=1,...,p-1$ 
and, thus,
$$
\mu(s,\la)\ \ge\ -a_1\alpha(p-1)p+a_2\alpha p.
$$
Now, the right hand expression is $>0$, by our choice of $a_1$ and $a_2$.
\par
Next, let $\eps=0$.
Since the definition of (semi)stability implies in that case 
$(\phi_{\id_{L^{\otimes r-1}}})\circ\cdots\circ\phi\neq 0$,
every one parameter subgroup acts with weight $\le 0$
on the ``${\bf P}$-component'' of $s$.
This allows us to argue in the same way as before.
\end{pf*}
\section{The moduli space of semistable Hitchin pairs}
\subsection{S-equivalence and the main result}
We define ${\cal M}^{(s)s}_{(L,P)}:={\frak P}^{(s)s}\catqot\SL(V)$.
Then ${\cal M}^{ss}_{(L,P)}$ is a projective scheme.
In order to describe its closed points, we have to introduce the
notion of S-equivalence:
For any semistable Hitchin pair, we can construct a Jordan-H\"older filtration
of $\E$
$$0=\E_0\subset\E_1\subset\cdots\subset\E_l=\E$$
by $\phi$-invariant subsheaves.
We obtain stable Hitchin pairs $(\E_i/\E_{i-1},\eps,\phi_i)$, $i=1,...,l$.
The associated graded object 
$$
\mathop{\rm gr}(\E,\eps,\phi):=\bigoplus (\E_i/\E_{i-1},\eps,\phi_i)
$$
is well-defined up to isomorphism.
We say that two semistable Hitchin pairs $(\E_1,\eps_1,\phi_1)$
and $(\E_2,\eps_2,\phi_2)$ of type $(L,P)$
are \it S-equivalent\rm , if the associated
graded objects are equivalent Hitchin pairs.
One can show that any semistable Hitchin pair degenerates into its
associated graded object and that the associated graded object is polystable.
We summarize the results of our discussions in:
\begin{Thm}
{\rm i)}
There is a natural transformation of functors
$$
\tau\colon M^{ss}_{(L,P)}\lra h_{{\cal M}_{(L,P)}^{ss}}
$$ such that for any other scheme $\tilde{\cal M}$ 
and any natural transformation
$\tau^\p\colon {\cal M}^{ss}_{(L,P)}\lra h_{\tilde{\cal M}}$
there is a uniquely determined morphism $\vartheta\colon {\cal M}^{ss}_{(L,P)}\lra \tilde{\cal M}$
with $\tau^\p=h(\vartheta)\circ \tau$.
\par
{\rm ii)}
${\cal M}^s_{(L,P)}$ is a coarse moduli space for the functor $M^s_{(L,P)}$.
\par
{\rm iii)}
The closed points of ${\cal M}^{ss}_{(L,P)}$ naturally correspond to the
S-equivalence classes of semistable Hitchin pairs of type $(L,P)$.
\end{Thm}
\subsection{The $\C^*$-action on ${\cal M}_{(L,P)}^{ss}$}
On the space ${\cal M}:={\cal M}_{(L,P)}^{ss}$ there is a natural
$\C^*$-action given by multiplication of $\phi$ by a constant.
The fixed point set is the union of the part which corresponds to the
Hitchin pairs $(\E,\eps,0)$, i.e., the Gieseker moduli space,
and the part ${\cal M}_\infty$ which corresponds to pairs $(\E,0,\phi)$.
The closed subset ${\cal M}_\infty$ is the part which compactifies 
the moduli space of semistable oriented Hitchin pairs.
Let ${\cal M}_{\neq 0}$ be the $\C^*$-invariant open subscheme
of semistable oriented Hitchin pairs, i.e., the set described by $\eps\neq 0$.
We observe that ${\cal M}_\infty={\cal M}_{\neq 0}\catqot \C^*$.
Here, we use that the GIT-quotient comes with a natural
ample line bundle  and that the $\C^*$ action is canonically linearized
in this line bundle.
\subsection{The Hitchin map}
Suppose that $X$ is a curve.
Let ${\frak P}^*$ be the open subset of the parameter space
${\frak P}$ parametrizing elements 
$([q\colon V\otimes \O_X\lra \E],[\eps,\phi])$ for which $\E$ 
is torsion free and $H^0(q)$ is an isomorphism,
and ${\frak P}^*_{\neq 0}$ the part of ${\frak P}^*$ lieing
in ${\frak Q}\times (V\otimes N)^\vee\otimes (V\otimes N\otimes M)$, 
i.e., the part parametrizing pairs
with $\eps\neq 0$.
Since the Quot scheme is reduced in this case, the restriction of 
${\frak E}_{\frak P}$ to ${\frak P}^*_{\neq 0}\times X$ is locally free.
This allows us to define the characteristical polynomial map
associated to $\phi_{{\frak P}|{{\frak P}^*_{\neq 0}\times X}}$:
$$
\chi_{{\frak P}^*_{\neq 0}}\colon {{\frak P}^*_{\neq 0}}\lra H^0(X,L^{\otimes r})
\oplus\cdots\oplus H^0(X,L).
$$
The $\C^*$-action on $(V\otimes N)^\vee\otimes (V\otimes N\otimes L)$
induces a $\C^*$-action on the right hand vector space which is given
on $H^0(X, L^{\otimes i})$ by multiplication with $z^i$, $i=1,...,r$.
Let $\C^*$ act on $\C$ by multiplication and form the weighted
projective space
$$
\widehat{\P}:=[H^0(X,L^{\otimes r})\oplus
\cdots\oplus H^0(X,L)\oplus\C]\catqot\C^*.
$$
Then the map $\chi_{{\frak P}^*_{\neq 0}}$ can be extended to a map
$
\chi_{\frak P}\colon {\frak P}^*\lra \widehat{\P}
$
which is invariant under the $\SL(V)$-action.
Thus we get a map
$$
\chi_{\M}\colon \M\lra \widehat{\P},
$$
which we call \it the Hitchin map\rm .
We oberve that $\chi_{\M}$ is proper by \cite{Ha}, II.4.8.(c),
applied to $f=\chi_{\M}$ and $g\colon \widehat{\P}\lra \{\text{pt}\}$.
\end{document}